\documentclass[12pt]{article}
\usepackage{amsfonts}
\usepackage{amsmath}
\usepackage{amssymb}
\usepackage{amsxtra}
\usepackage{graphicx}
\usepackage{geometry}
\usepackage{color}
\usepackage{colortbl}
\newcommand{\RR}{{\rm I\hskip -.15em R}}
\def\rd{\RR^d}
\newfont{\ruk}{eusm10 at 12pt}

\def\div{{\rm div}}
\def\lm{\lim\limits_{\ve\to 0}}
\def\O{\Omega}
\def\beq{\begin{equation}}
\def\eeq{\end{equation}}

\def\nab{\nabla}
\def\f{\frac}

\def\ve{\varepsilon}
\def\vp{\varphi}
\def\dx{\,dx}

\def\ild{\int_{\rd}}

\def\ilb{\int_{Y}}

\def\ds{\displaystyle}

\def\l{\left}
\def\r{\right}
\def\ue{u^\ve}

\def\C0{C_0^\infty}

\def\pa{\partial}

\newtheorem{theorem}{Theorem}[section]
\newtheorem{lemma}{Lemma}
\newtheorem{opr}[theorem]{Definition}

\newcommand{\doc}{{\bf Proof. }}

\def\a{\alpha}
\def\b{\beta}

\chardef\No=194

\def\dx{\,dx}

\def\per{{\rm per}}

\def\la{\lambda}

\def\*{^*}

\def\n{\noindent}
\def\E{\mathcal{E}}
\def\e{\varepsilon}
\begin{document}
\begin{center} 
\begin{Large} 
\begin{bf}

      Estimates in  homogenization  \\[0.2cm]
       of higher-order elliptic operators\\[.7cm]
       
\end{bf} 
\end{Large}

\begin{large}
    
        Svetlana E. PASTUKHOVA\\[0.5cm]                        
        {\em MIREA (Technical University)}\\
        {\em Department of Mathematics}\\
        {\em prospect Vernadskogo 78}\\
        {\em 119454 Moscow, Russia}\\[0.2cm]
       
\end{large}
\end{center}
{\em Abstract.}
\begin{small} 
A divergent-type elliptic operator $A^\e$ of arbitrary even order $2m$ is studied. Coefficients of the operator are $\e$-periodic, $\e{>}0$ is a small parameter. The resolvent equation $A^\e u^\e{+}\lambda u^\e{=}f$ is solvable in the Sobolev space $H^m(\rd)$ of order $m$ for any $f{\in}L^2(\rd)$, provided  the parameter $\la$ is sufficienly large, $\la{>}\Lambda$, where the bound $\Lambda$ depends only on constants from ellipticity condition. The limit 
equation is of the same type but with constant coefficients, that is, $\hat{A}  u{+}\lambda u{=}f$. The limit operator $\hat{A}$ can be considered here, for instance, in the sense of ${G}$-convergence. We prove that the resolvent $(\hat{A}{+}\lambda)^{-1}$ approximates $(A^\e{+}\lambda)^{-1}$ in operator $(L^2{\to}L^2)$-norm with the estimate 
$\|(A^\e{+}\lambda)^{-1}-(\hat{A}{+}\lambda)^{-1}\|_{L^2(\rd){\to}L^2(\rd)}=O(\e)$, as $\e{\to}0$. We find also the approximation of 
the resolvent $(A^\e{+}\lambda)^{-1}$ in operator $(L^2{\to}H^m)$-norm. This is the sum 
$(\hat{A}{+}\lambda)^{-1}{+}\mathcal{K}^\e$, where $\mathcal{K}^\e$ is a correcting operator whose structure is given.
We prove the estimate $\|(A^\e{+}\lambda)^{-1}-(\hat{A}{+}\lambda)^{-1}{-}\mathcal{K}^\e\|
_{L^2(\rd){\to}H^m(\rd)}=O(\e)$, as $\e{\to}0$.
\end{small}
\section{Introduction}

The theory of $G$-convergence of differential operators and connected with it the theory of
multidimensional homogenization have being studied from long ago since 60-th, see e.g. \cite{Sp}, \cite{GSp}.
Firstly,
there was investigated the class $\E(\lambda_0,\lambda_1)$ of second order elliptic operators of the form
$$A{=}\sum_{i,j=1}^d D_i(a_{ij}D_j)\,\, (a_{ij}{=}a_{ji}),$$ where 
$a_{ij}(x)$ are measurable functions in a bounded domain $\O\subset\rd$ subject to the inequality
\[
\lambda_0 |\xi|^2\le \sum_{i,j=1}^d a_{ij}(x)\xi_i\xi_j\le
\lambda_1 |\xi|^2\quad (x\in\O, \xi\in \rd),
\]
$\lambda_0$, $\lambda_1 $ are positive constants. The sequence of operators $A_k{\in}\E(\lambda_0,\lambda_1),\, k{\in}\mathbb{N},$ is called 
$G$-convergent to the operator $A{\in}\E(\lambda_0,\lambda_1)$, as $k{\to} \infty$, if for any $f{\in}L^2(\O)$ the sequence of solutions $u_k$ of the Dirichlet problem 
$A_ku_k{=}f\,\, (u_k{\in}H_0^1(\O))$ converges in $L^2(\O)$ to the 
solution $u$ of the Dirichlet problem 
$Au{=}f\,\, (u{\in}H_0^1(\O))$. In \cite{Sp} there was proved that the class  $\E(\lambda_0,\lambda_1)$  is compact in the sense of $G$-convergence.

In  \cite{ZKOK} there were introduced classes of divergent-type elliptic operators of arbitrary order $2m{\ge}2$
for which the
 $G$-compactness was proved. Numerous properties of $G$-convergence in these classes were established, among them  the
construction receipt for $G$-limiting operator  of the family of higher-order differential operators having $\e$-periodic coefficients when the small parameter $\e$ is tending to zero.
The latter result relates to the homogenization theory and may be amplified further.
Now, we are interested in the rate of convergence 
of the solution to initial nonhomogeneous problem to the 
solution of homogenized problem 
going with all this,   for the time being, not in a bounded domain but
in the whole space $\rd$.
We prove that the rate of this convergence is of order $\e$. We assume minimal regularity conditions on the data of the problem. 
Thereby, the result may acquire the form of operator-type convergence of resolvents in terms of operator norms with corresponding estimate on the rate of convergence.
This kind estimate
 for operators in the class $\E(\lambda_0,\lambda_1)$ defined above is known and was first proved in \cite{BS}.
 To obtain the 
 estimate in homogenization of higher-order differential operators, we use 
 the approach proposed in \cite{Z05}, \cite{ZP}, certainly, necessarily modified.  
 Formerly, this approach was applied only to second order differential equations, though of different types, among them, equations with various properties of degeneracy, with multi-scale coefficients or quasiperiodic coefficients,  non-linear and non-divergent equations, system of elasticity theory equations and parabolic equations \cite{ZPpar}, \cite{SMZ},  \cite{PT}, \cite{MIAN}, \cite{FA}, \cite{RJMPh15}, \cite{CaPZ}, \cite{ET}.
 
 As one of
 distinctive features of our problem we mention also  a presence of lower terms 
  combined with a  lack of any 
 symmetry in coefficients matrix.
  This allows to take along divergent-type higher-order differential operators of general form. Maybe, for even  second order elliptic operators with non-symmetric  coefficients matrix with lower terms, the obtained result is interesting.
 
 Certainly, for applications it is important, above all, to prove 
 operator-type estimates in homogenization of boundary problems in bounded domains. According to the method for derivation of such estimates, given in   \cite{Z05}, \cite{ZP}, one should 
 start from resolvent equations in the whole space. Approximations constructed at the first stage to satisfy merely the equation 
 are then adjusted to  the boundary conditions by means of additional correctors which are of boundary layer nature.
 So, this paper may be regarded as a preparatory one for future derivation of homogenization estimates in problems with 
  boundary conditions which are of the main interest for us.

\section{  Main results}

\n \textbf{2.1.}
We begin with  some preliminary material. 

Denote by $H^m{=}H^m(\rd)$, $m{\ge}0$ is integer,
the Sobolev space equipped with a norm defined by the equality 
\[
\|u\|^2_{H^m}=\ild \sum\limits_{|\alpha|\le m}|D^{\alpha}u|^2\dx,
\]
where $\alpha{=}(\alpha_1,\ldots,\alpha_d)$
is 
a  multi-index with  non-negative integer components $\alpha_i$,
\[
|\alpha|{=}\sum_{i=1}^d\alpha_i,\quad D^{\alpha}=D_1^{\alpha_1}\ldots D_d^{\alpha_d},\quad D_i=\f{\partial}{\partial x_i},
\]
$H^0(\rd){=}L^{2}(\rd)$. It is known that  the set $C_0^\infty(\rd)$ is dense in $H^m(\rd)$. Let
\[
\|u\|_j=(\ild \sum\limits_{|\alpha|= j}|D^{\alpha}u|^2\dx)^{1/2} 
\]
for integer $j{\ge}0$, then the expression $(\|u\|_m^2{+}\|u\|_0^2)^{1/2}$ defines the equivalent norm in $H^m(\rd)$, $m{\ge}1$.

Denote by $H^{-m}{=}H^{-m}(\rd)$, $m{\in}\mathbb{N}$, the space dual to $H^m(\rd)$. We have a triple of spaces
\beq\label{0}
H^m(\rd)\subset L^{2}(\rd)\subset H^{-m}(\rd) 
\eeq 
for each integer $ m>0.$ 

Throughout this paper, $C_0^\infty(\rd)$, 
$ L^{2}(\rd)$, $H^m(\rd)$ are considered as spaces of real-valued functions.

Consider linear differential operators of the form
\beq\label{1}
A {=}
\sum\limits_{|\alpha|,|\beta|\le m}(-1)^{|\a|}D^{\alpha}(a_{\alpha\beta}(x)D^{\beta})
\eeq 
with bounded measurable real-valued coefficients $a_{\alpha\beta}(x)$. 
\begin{opr}\label{Opr1}
 A  matrix $\{a_{\alpha\beta}(x)\}_{|\alpha|,|\beta|\le m}$ is called 
 elliptic 
 if
\beq\label{2}\|a_{\alpha\beta}\|_{L^\infty(\rd)}\le \lambda_1,\quad |\alpha|,|\beta|\le m,
\eeq
and
\beq\label{3}
\int_{\rd}\sum\limits_{|\alpha|=|\beta|=m}a_{\alpha\beta}(x)D^{\beta}u
 D^{\alpha}u\,dx\ge  \lambda_0
\int_{\rd} \sum\limits_{|\alpha|=m}| D^{\alpha}u|^2\,dx\quad \forall u\in C_0^\infty(\rd)
\eeq
for some constants $\lambda_1,\,\lambda_0>0$.
\end{opr}

It is rather known that the necessary condition for (\ref{3}) is the following algebraic inequality
\beq\label{03}
\sum\limits_{|\alpha|=|\beta|=m}a_{\alpha\beta}(x)\xi^{\beta}
 \xi^{\alpha}\ge  \lambda_0
 \sum\limits_{|\alpha|=m}| \xi^{\alpha}|^2\quad \forall \xi\in \rd,\quad \xi^\a=\xi_1^{\a_1}\ldots\xi_d^{\a_d},
\eeq
which holds for a.e. $x{\in}\rd$. 
In the case of the constant matrix $\{{a}_{\a\b}\}$, the coerciveness inequality of the form (\ref{3}) is equivalent to (\ref{03}). The latter can be readily shown  due to the  Plancherel identity 
via Fourier transform images.

\bigskip
The differential expression of the form (\ref{1}), corresponding to an elliptic matrix $\{a_{\alpha\beta}(x)\}$, defines a bounded linear operator $A: H^{m}\to H^{-m}$ in the following way.
Given $u{\in}H^m$ and $f{\in}H^{-m}$, we say that $Au{=}f$ if the action of $f$
is described  as 
\[
\langle f,\vp\rangle=
\sum\limits_{|\alpha|\le m,|\beta|\le m}(a_{\alpha\beta}D^{\beta}u,D^{\alpha}\vp),\quad \vp{\in}\C0(\rd),
\]
where 
$(u,\vp)=
(u,\vp)_{L^{2}(\rd)}$ and $\langle f,\vp\rangle$  being the value of  $f{\in}H^{-m}$ on the element $\vp{\in}H^m$.
Clearly, $\|A\|\le \lambda_1.$ Moreover, the operator $A: H^{m}{\to} H^{-m}$ is lower semibounded  and 
satisfies the inequality 
\beq\label{4}
\langle Au,u\rangle\ge \lambda_0/2\|u\|^2_m
- \lambda_2\|u\|^2_0 
 \quad \forall u\in C_0^\infty(\rd),
\eeq
where $\lambda_2{\ge}0$ is a  constant depending on $\lambda_0$ and $\lambda_1$ from the property ellipticity.
In fact,
\[
\ds {\langle Au,u\rangle=\sum\limits_{|\alpha|=|\beta|=m}(a_{\alpha\beta}D^{\beta}u,D^{\alpha}u)+
\sum\limits_{|\alpha|+|\beta|<2m}(a_{\alpha\beta}D^{\beta}u,D^{\alpha}u)\ge
}
\atop\ds{\lambda_0\|u\|^2_m-\lambda_1(\delta\|u\|^2_m+C_\delta\sum_{j=0}^{m-1}\|u\|^2_j) \quad \forall\,\delta>0,
}
\]
where in its turn
\[
\|u\|^2_j\le \tau
\|u\|^2_m+C_\tau \|u\|^2_0,\quad 1\le j\le m-1,\quad \forall\,\tau{>}0.
\]
Hence, (\ref{4}) is  easily obtained if  $\delta$ and $\tau$ are chosen sufficiently small. 
From (\ref{4}), we derive that the operator $A{+}\lambda I: H^{m}\to H^{-m}$, $\lambda{\ge}\lambda_2{+}1$, is coercive and satisfies the inequality
\[
\langle( A{+}\lambda I)u,u\rangle\ge \lambda_0/2\|u\|^2_m+\|u\|^2_0, 
\]
for any $u\in C_0^\infty(\rd)$. In terms of
the standard norm in $H^m$   it means 
\[
\langle( A{+}\lambda I)u,u\rangle\ge \tilde{\lambda}_0\|u\|^2_{H^m}, \quad 
\tilde{\lambda}_0=const(\lambda_0).
\]

Consequently, we come to the following conclusion.
\begin{lemma}\label{L1}
Let $A$ be an operator of form (\ref{1}) where a coefficients matrix is elliptic with   constants $\lambda_1,\,\lambda_0$.
Then for any $f{\in}  H^{-m}$ and sufficiently large positive $\lambda$,
the equation
\beq\label{5}
( A{+}\lambda I)u=f, 
\eeq
has the unique solution $u{\in}  H^{m}$, and
\beq\label{6} \|u\|_{H^m}= \|( A{+}\lambda I)^{-1}f\|_{H^m}\le
C\|f\|_{H^{-m}},
\quad 
C=const(\lambda_1,\lambda_0).\eeq 
In other words, the elliptic differential operator (\ref{1}) 
determines an isomorphism $( A{+}\lambda I){:}H^m{\to} H^{-m}$ for each $\lambda{\ge}\Lambda$,  and any functional $f{\in}  H^{-m}$ admits the unique representation 
$f{=}( A{+}\lambda I)u$, where $u{\in}  H^{m}$.
\end{lemma}

This  statement is based on the following result for abstract operators (see \cite{ZKOK}, Theorem 1 in Chapter  I) which is often called Lax--Milgram Theorem.

Let $V$ be a real reflexive 
Banach space and $V^\prime$ its dual. The value of a functional $f{\in}V^\prime$ at a point $v{\in}V$ is denoted by  $\langle f,v\rangle$.

\noindent\textbf{Proposition.} \textit{Let } $L:V\to V^\prime$
\textit{be a continuous linear operator such that } $\langle Lv,v\rangle\ge \tilde{\lambda}\|v\|^2_V$
\textit{for any } $v{\in}V$.
\textit{Then the equation } $Lv{=}f\,\, (v{\in}V)$
\textit{is uniquely solvable for any } $f{\in}V^\prime$,
\textit{and } 
\[
 \|v\|_V=\ \|L^{-1}f\|_{V}\le\tilde{\lambda}^{-1}\|f\|_{V^\prime}.
\]

\bigskip
In accordance  with (\ref{6}) and due to the embedding (\ref{0}), the resolvent $( A{+}\lambda I)^{-1}$ 
can be considered as an operator in $L^2(\rd)$. Evidently, 
$$\|( A{+}\lambda I)^{-1}\|_{L^2(\rd)\to L^2(\rd)}\le \tilde{\lambda}_0^{-1}, \quad 
\lambda{\ge}\lambda_2.$$

\bigskip
\n \textbf{2.2.} Consider a family of differential operators depending on a small parameter $\e{>}0$ 
\beq\label{8}
\ds{
A^\e =
\sum\limits_{|\alpha|,|\beta|\le m}(-1)^{|\a|}D^{\alpha}(a^\e_{\alpha\beta}(x)D^{\beta}),}
\atop\ds{ a^\e_{\alpha\beta}(x)=a_{\alpha\beta}(x/\e),}
\eeq 
the corresponding matrix $\{a_{\alpha\beta}(x)\}$ is elliptic (in the sense of Definition \ref{Opr1}) and  is 1-periodic in each variable $y_1,\ldots,y_d$, $Y{=}[-\f12,\f12)^d$ is the periodicity cell.

It is clear that the matrix $\{a^\e_{\alpha\beta}(x)\}$ is elliptic (in the sense of Definition \ref{Opr1}) for any  $\e{\in}(0,1).$
In particular, the uniform in $\e$ coerciveness estimate holds
\beq\label{pe1}
\int_{\rd}\sum\limits_{|\alpha|=|\beta|=m}a^\e_{\alpha\beta}(x)D^{\beta}u
 D^{\alpha}u\,dx\ge  \lambda_0
\int_{\rd} \sum\limits_{|\alpha|=m}| D^{\alpha}u|^2\,dx\quad \forall u\in C_0^\infty(\rd).
\eeq
By closure, this estimate is true for any $u\in H^m(\rd).$

According to the theory of $G$-convergence of elliptic operators, developed in  \cite{ZKOK},
there is a strong $G$-convergence of $A^\e$ to the limit operator $\hat{A}$,
\beq\label{9}
A^\e\stackrel{G}\Longrightarrow \hat{A}, \quad 
\hat{A} =
\sum\limits_{|\alpha|,|\beta|\le m}(-1)^{|\a|}D^{\alpha}(\hat{a}_{\alpha\beta}D^{\beta}),
\eeq 
$\hat{A}$ is  of the form (\ref{1}) with a constant elliptic matrix $\{\hat{a}_{\alpha\beta}\}$. 
The way how to find  the  limit matrix $\{\hat{a}_{\alpha\beta}\}$ is described later in  Sect. 3 (see (\ref{c4}),(\ref{c1})).
For the definition of the strong $G$-convergence, used  in (\ref{9}),   and also for its properties see \cite{ZKOK}.
We shall extend this result  in another direction, or maybe, look at the limit operator from another side, in the sense of somehow stronger convergence.

Due to the above arguments (see Lemma \ref{L1}), since  both operators $A^\e$ and $\hat{A} $
are elliptic, the resolvents $(A^\e+\lambda I)^{-1}$ and $(\hat{A}+\lambda I)^{-1} $
exist for sufficiently large  $\lambda$, $\lambda{\ge}\Lambda(\lambda_0,\lambda_1)$, where $\lambda_0,\lambda_1$ are the constants from the ellipticity condition. As operators in $ L^2(\rd)$, these resolvents turn to be close to each other in the operator norm, and the degree of closeness is of order $\e$.
\begin{theorem} Under the above assumption of ellipticity there holds the estimate
\beq\label{10}
\|(A^\e+\lambda I)^{-1}-(\hat{A}+\lambda I)^{-1}\|_{L^2(\rd)\to L^2(\rd)}\le c_0 \e , \quad c_0=const(\lambda_0,\lambda_1),\quad \lambda{\ge}\Lambda(\lambda_0,\lambda_1).
\eeq
\end{theorem}

 In \cite{V} the estimate (\ref{10}) was proved for self-adjoint operators (\ref{8}) without lower order terms. To this end the spectral approach from \cite{BS} was used. 

To give a simple interpretation of (\ref{10}) consider 
equations
\beq\label{11}
u^\e\in H^m(\rd),\quad A^\e u^\e+\lambda  u^\e=f, 
\eeq
\beq\label{12}
u\in H^m(\rd),\quad \hat{A}u+\lambda u=f, 
\eeq
for an arbitrary $f{\in}L^2(\rd)$. The operator estimate (\ref{10}) means that 
\beq\label{13}
\|u^\e-u\|_{L^2(\rd)}\le \e c_0  \|f\|_{L^2(\rd)} , \quad c_0=const(\lambda_0,\lambda_1),
\eeq
provided $\la$ is sufficiently large.

One can consider the resolvent $(A^\e+\lambda I)^{-1}$ as operator $ L^2(\rd){\to} H^m(\rd)$. Then its approximation in the operator  $( L^2{\to} H^m)$-norm should be taken as a sum of 
 $(\hat{A}+\lambda I)^{-1}$ and the correcting operator $\mathcal{K}^\e$ whose structure is easily restored through the formula (\ref{2.36}).
Moreover, there holds the estimate
\beq\label{14}
\|(A^\e+\lambda I)^{-1}-(\hat{A}+\lambda I)^{-1}-\mathcal{K}^\e \|_{L^2(\rd)\to H^m(\rd)}\le c_0 \e , \quad c_0=const(\lambda_0,\lambda_1). 
\eeq
Certainly,
  the constant $c_0$ in (\ref{10}) and (\ref{14}) depends also on the dimension $d$
  and the order $m$, but this will not be mentioned  anymore. 
 
  The direct proof of the above estimates   is given in Sections 5 and 6. Some necessary preliminaries are carried away to Sections 3 and 4.


\section{ Cell problems}
\setcounter{equation}{0}
  \setcounter{theorem}{0}
   \setcounter{lemma}{0}
  
 \n \textbf{3.1.}
  On the set of periodic  functions $u{\in}C_\per^\infty(Y)$ with zero mean, that is 
\[
\langle u\rangle=\ilb u(y)\,dy=0,
\]  
the expression 
\[
\left(\ilb \sum_{|\alpha|=m}|D^{\alpha}u|^2\,dy\right)^{1/2}
\]
yields   the norm. We denote  by 
$\mathcal{W}$  the completion of this set under the norm in question.

The estimate (\ref{3}) for the  periodic matrix $\{a_{\alpha\beta}(x)\}$ leads to the following inequality for 
periodic functions
\beq\label{c2}
\ilb\sum\limits_{|\alpha|=|\beta|=m}a_{\alpha\beta}(y)D^{\beta}u
 D^{\alpha}u\,dy\ge  \lambda_0
\ilb \sum\limits_{|\alpha|=m}| D^{\alpha}u|^2\,dy\quad \forall u\in C_\per^\infty(Y),
\eeq
which is verified a little bit later.
By closure, (\ref{c2}) holds for any  $u{\in}\mathcal{W}$ and means that
\[
A^0=\sum_{|\alpha|=|\beta|=m}D^{\alpha}(a_{\alpha\beta}(y)D^{\beta})
\] 
is a coercive operator $\mathcal{W}\to \mathcal{W}^\prime$.

For each multi-index $\gamma$ with $|\gamma|{\le}m$ consider the equation
\beq\label{c1}
N_\gamma\in \mathcal{W}, \quad
\sum_{|\alpha|=|\beta|=m}D^{\alpha}(a_{\alpha\beta}(y)D^{\beta}N_\gamma(y))=-\sum_{|\alpha|=m}D^{\alpha}(a_{\alpha\gamma}(y)),
\eeq  
where $\{a_{\alpha\beta}(y)\}$ is the periodic matrix from (\ref{8}).
 The right-hand side in (\ref{c1}) determines in a natural way a linear functional  $F_\gamma$ on $\mathcal{W}$, and the equation itself can be written in the operator form
$ 
A^0N_\gamma{=}F_\gamma\, (N_\gamma{\in} \mathcal{W}).$
 Therefore, the unique solvability of (\ref{c1})  follows from  Proposition  given after Lemma \ref{L1}.
 
 Return to the assertion which is  basic for this Section.
  \begin{lemma}\label{lem3.2} Suppose that 
 $\{a_{\alpha\beta}(x)\}$ is an elliptic periodic  matrix.  Then
the estimate (\ref{3}) entails the 
 inequality (\ref{c2}).
\end{lemma}
\doc
Substitute
in (\ref{pe1})  the finite function
\[
\psi_\e(x)=\e^m v(x/\e)\vp(x),\quad v\in C^\infty_\per(Y),\quad \vp\in \C0(\rd),
\] 
such that
\beq\label{pe2}
D^\a\psi_\e(x)=(\pa^\a v)^\e(x)\vp(x)+O(\e),
\quad (\pa^\a v)^\e(x)=(D^\a_y v(y))|_{y=x/\e},\quad |\a|=m.
\eeq
Here and hereafter,  $\pa^\a{=}D^\a_y$.  
Obviously,  we obtain 
\[
\int_{\rd}\sum\limits_{|\alpha|=|\beta|=m}a^\e_{\alpha\beta}(x)(\pa^\b v)^\e (x)(\pa^\a v)^\e(x) |\vp(x)|^2
\,dx{\ge}  \lambda_0
\int_{\rd} \sum\limits_{|\alpha|=m} |(\pa^\a v)^\e(x)|^2 |\vp(x)|^2\,dx{+}O(\e),
\]
that after the passage to the limit, as $\e\to 0 $,  gives
\[
\langle\sum\limits_{|\alpha|=|\beta|=m} a_{\alpha\beta}\pa^\b v \pa^\a v\rangle
\int_{\rd} |\vp(x)|^2
\,dx{\ge}  \lambda_0 \langle\sum\limits_{|\alpha|=m} |\pa^\a v|^2\rangle 
\int_{\rd}  |\vp(x)|^2\,dx
\]
and, finally,
\[
\langle\sum\limits_{|\alpha|=|\beta|=m} a_{\alpha\beta}\pa^\b v \pa^\a v\rangle
\ge  \lambda_0
 \langle\sum\limits_{|\alpha|=m} |\pa^\a v|^2\rangle 
\]
which is exactly the inequality (\ref{c2}).
 Above, the so-called mean value property of periodic functions is used: 

\textit{ if 
$b{\in} L^1_\per(Y)$, $b^\e(x){=}b(y)|_{y=x/\e}$, $\psi{\in}\C0(\rd),$ then 
}
\beq\label{mvp}
\lm \int_{\rd}b^\e(x)\psi(x)\,dx=\langle b\rangle \int_{\rd}\psi(x)\,dx,\quad \langle b\rangle = \ilb b(y)\,dy.
\eeq

\bigskip
\n \textbf{3.2.}
The coefficients  of the operator $\hat{A}$
(see (\ref{9})) are defined 
with the help of the solutions of 
cell problems (\ref{c1}),
\beq\label{c4}
\hat{a}_{\alpha\beta}=
\langle a_{\alpha\beta}(y)+\sum\limits_{|\gamma|=m}
a_{\alpha\gamma}(y)D^{\gamma}N_\beta(y)\rangle,\quad |\alpha|\le m,|\beta|\le m.
\eeq
Introducing the symbol $e_{\alpha\beta}$
with multi-indices $\alpha$, $\beta$
such that
$$
e_{\alpha\beta}
=\l\{
\begin{array}{rcl}
1,\text{ if }\alpha=\beta,\\
0\text{ otherwise },\\
\end{array}\r.
$$
we rewrite
\beq\label{c5}
\hat{a}_{\alpha\beta}=
\langle \sum\limits_{|\gamma|=m}
a_{\alpha\gamma}(y)(e_{\beta\gamma}+D^{\gamma}N_\beta(y))\rangle, 
\eeq 
or 
\beq\label{c6}
\hat{a}_{\alpha\beta}=\langle \tilde{a}_{\alpha\beta}\rangle,\quad \tilde{a}_{\alpha\beta}(y)=
 \sum\limits_{|\gamma|=m}
a_{\alpha\gamma}(y)(e_{\beta\gamma}+D^{\gamma}N_\beta(y)).
\eeq

\begin{lemma}\label{lem3.1} The matrix $\{\hat{a}_{\a\b}\}$, defined by relations (\ref{c4}), (\ref{c1}), is elliptic.
\end{lemma}
\doc 
Evidently, we need to verify for  $\{\hat{a}_{\a\b}\}$ only the coerciveness property from the definition of elliptic matrices, that is,
\beq\label{cp}
\int_{\rd}\sum\limits_{|\alpha|=|\beta|=m}\hat{a}_{\alpha\beta}D^{\beta}w
 D^{\alpha}w\,dx\ge  \lambda_0
\int_{\rd} \sum\limits_{|\alpha|=m}| D^{\alpha}w|^2\,dx\quad \forall w\in C_0^\infty(\rd).
\eeq
To this end, substituting in (\ref{pe1}) the function
\beq\label{cp1}\ds{
u_\e(x)=w(x)+\e^m\sum_{|\gamma|=m}N^\e_\gamma(x)D^\gamma w(x),}
\atop\ds{w\in \C0(\rd),\quad N^\e_\gamma(x)=N_\gamma(y)|_{y=x/\e},
}
\eeq
where $N_\gamma(y)$ is the solution of  the cell problem (\ref{c1}), we obtain
\beq\label{pe100}
\int_{\rd}\sum\limits_{|\alpha|=|\beta|=m}a^\e_{\alpha\beta}(x)D^{\beta}u_\e
 D^{\alpha}u_\e\,dx\ge  \lambda_0
\int_{\rd} \sum\limits_{|\alpha|=m}| D^{\alpha}u_\e|^2\,dx\quad \forall w\in C_0^\infty(\rd)
\eeq
and pass here to the limit, as $\e{\to}0$.

First calculate 
\[\ds{ \sum_{|\b|=m}{a}^\e_{\a\b}D^\b u_\e
\stackrel{(\ref{cp1})}=\sum_{|\b|=m}{a}^\e_{\a\b}(D^\b w+
\sum_{|\gamma|=m}(\pa^\b N _\gamma)^\e D^\gamma w) +r^\e_\a=
}
\atop\ds{\sum_{|\b|=m}\sum_{|\gamma|=m}{a}^\e_{\a\gamma}(e_{\b\gamma}+(\pa^\gamma N _\b)^\e)D^\b w+r^\e_\a,
}
\]
\[\ds{ D^\a u_\e
\stackrel{(\ref{cp1})}=D^\a w+
\sum_{|\delta|=m}(\pa^\a N _\delta)^\e D^\delta w) +\tilde{r}^\e_\a=
}
\atop\ds{\sum_{|\delta|=m}(e_{\a\delta}+(\pa^\a N _\delta)^\e)D^\delta w+\tilde{r}^\e_\a,\quad |\a|=m,
}
\]
where the notation from (\ref{pe2}) is used and ${r}^\e_\a$, $\tilde{r}^\e_\a$ denote terms of order $O(\e)$.
Hence, 
\[\ds{ \sum_{|\a|=|\b|=m}{a}^\e_{\a\b}D^\b u_\e D^\a u_\e=
}
\atop\ds{\sum_{|\b|=|\delta|=m}[\sum_{|\a|=|\gamma|=m}{a}^\e_{\a\gamma}(e_{\b\gamma}+(\pa^\gamma N _\b)^\e)
(e_{\a\delta}+(\pa^\a N _\delta)^\e)]
D^\b w D^\delta w+R^\e,\quad R^\e=O(\e).
}
\]
By the  mean value property of periodic functions (see (\ref{mvp}))
\beq\label{cp2}
\ds{ \lm\int_{\rd}\sum_{|\a|=|\b|=m}{a}^\e_{\a\b}D^\b u_\e D^\a u_\e\dx=
}
\atop\ds{\int_{\rd}\sum_{|\b|=|\delta|=m}\langle\sum_{|\a|=|\gamma|=m}{a}_{\a\gamma}(e_{\b\gamma}+D^\gamma N _\b)
(e_{\a\delta}+D^\a N _\delta)\rangle
D^\b w D^\delta w\dx=
}
\eeq
\[
\int_{\rd}\sum_{|\b|=|\delta|=m}\hat{a}_{\a\b}
D^\b w D^\delta w\dx,
\]
where  the mean over $Y$ in (\ref{cp2}) is  reduced to $\hat{a}_{\a\b}$ by using the definition of homogenized coefficients and the equality
\[
\langle\sum_{|\a|=|\gamma|=m}{a}_{\a\gamma}(e_{\b\gamma}+D^\gamma N _\b)
D^\a N _\delta)\rangle=0,
\]
which stems from (\ref{c1}).

Again with the help of the  mean value property, we can show the weak convergence 
 $u_\e\rightharpoonup w$ in $H^m(\rd)$. Here,  the structure of $u_\e$ is essential (see (\ref{cp1})). So, there is the lower semi-continuity property of the norm
\beq\label{cp3}
\liminf_{\e\to 0}\int_{\rd} \sum\limits_{|\alpha|=m}| D^{\alpha}u_\e|^2\,dx\ge
\int_{\rd} \sum\limits_{|\alpha|=m}| D^{\alpha}w|^2\,dx.
\eeq 

From (\ref{pe100})--(\ref{cp3}) it follows that (\ref{cp}) is true. The lemma is proved.

Since $\hat{A}$  is elliptic and with constant coefficients, we conclude that the solution $u$ of the homogenized equation 
(\ref{12})  really belongs to $H^{2m}(\rd)$ and the following estimate holds
\beq\label{ele}
\|u\|_{H^{2m}(\rd)}\le  C \|f\|_{L^{2}(\rd)}, \quad C=const(\la_0,\la_1).
\eeq

\section{On representation of solenoidal vectors}
\setcounter{equation}{0}
  \setcounter{theorem}{0}
  \setcounter{lemma}{0}

It is known that periodic solenoidal vectors with zero mean admit representation via  
divergence of a periodic skew-symmetric matrix. More precisely, \textit{ for any vector $g{\in} L^2_\per(Y)$, such that $\div\, g{=}0$ and $\langle g\rangle{=}0$,  there exists a skew symmetric 
matrix $G{\in} H^1_\per(Y)$ such that $\div \, G{=}g$ and } $\|G\|_{H^1_\per(Y)}\le c\|g\|_{H^1_\per(Y)}$, $c{=}const(d)$ (see the proof in \cite{JKO}, Chapter I,\S1).
Here, the relation 
\beq\label{p0}
\div\, g=\Sigma_{i=1}^dD_ig_i=0\eeq
for the vector $g{=}\{g_i\}{\in} L^2_\per(Y)^d$ means that
\beq\label{p1}
\langle g\cdot\nabla \varphi\rangle=0\quad\forall \,\varphi\in C_\per^\infty(Y).
\eeq
 The following lemma extends the above assertion to the situation when instead of operator of gradient $\nabla{:}H^1_\per(Y){\to}  L^2_\per(Y)^d$ and its 
 adjoint operator (that is the operator of divergence $\div$),
 one consider a pair of their analogues adequate to differential equations of order $2m$. That is, first,
 the operator of gradient of order $m$, acting from $ H^m_\per(Y)$ to vector-valued space $ L^2_\per(Y)^p$,
 $p$ being the number of multi-indices $\a$ with $|\a|{=}m$, and, second, its adjoint operator which may be called  "divergence of order $m$"${}$. Certainly, there appear "solenoidal vectors"${}$ corresponding to this type of divergence (see  below (\ref{lp1})$_2$), for which analogues of (\ref{p0}), (\ref{p1})) are fulfilled.

\begin{lemma}\label{lem4.1} Let
$\{g_\alpha\}_{|\alpha|=m}$ be a 1-periodic vector from $L^2(Y)^p$ 
such that
\beq\label{lp1}
\langle g_\alpha\rangle=0,\quad
\sum_{|\alpha|=m}D^\alpha g_\alpha=0.
\eeq
Then there is a 1-periodic matrix $\{G_{\alpha\beta}\}_{|\alpha|=|\beta|=m}$ from $H^m(Y)^{p\times p}$
such that for any $\alpha$, $\beta$
\beq\label{lp2}
G_{\alpha\beta}=-G_{\beta\alpha},
\eeq
\beq\label{lp3}
\|G_{\alpha\beta}\|_{H^m(Y)}\le c \sum_{|\alpha|=m}\|g_\alpha\|_{L^2(Y)},\quad c=const( {d,m}),
\eeq
\beq\label{lp4}
\sum_{|\gamma|=m}D^\gamma G_{\alpha\gamma}=g_\alpha.
\eeq
\end{lemma}
\doc Each component $g_\alpha$ admits  Fourier decomposition
\[
g_{\alpha}(y)=\sum_{0\neq n\in \mathbb{Z}^d}g^{n}_\alpha e^{2\pi i n\cdot y}, \quad  i=\sqrt{-1},
\]
and   there holds Parseval's identity
\[
\|g_{\alpha}\|^2_{L^2(Y)}=\sum_{n\in \mathbb{Z}^d}|g^{n}_\alpha|^2.
\]
Withal the equation (\ref{lp1})$_2$ implies that
\beq\label{lp5}
\sum_{|\alpha|=m}n^\alpha g^n_\alpha=0,\quad 0\neq n\in \mathbb{Z}^d.
\eeq

Define $G_{\alpha\beta}$ 
 via Fourier decomposition with coefficients
\[\ds{
G_{\alpha\beta}^n=(-g_{\beta}^nn^\alpha+g_{\alpha}^nn^\beta)\Lambda_m(n)^{-1}
(2\pi i)^{-m},\quad 0\neq n\in \mathbb{Z}^d,}
\atop\ds{\Lambda_m(n)=\sum_{|\gamma|=m}n^\gamma n^\gamma.
}
\]
By construction,
condition (\ref{lp2}) is evidently fulfilled 
and
\[
\sum_n |G_{\alpha\beta}^n n^\gamma|^2\le C \sum_{n\in \mathbb{Z}^d}|g^{n}_\alpha|^2,\quad |\gamma|=m,
\] 
thereby, the estimate (\ref{lp3}) holds true.
Moreover, for each $ 0\neq n\in \mathbb{Z}^d$
\[\ds{
e^{-2\pi i n\cdot y}
\sum_{|\beta|=m}D^\beta[G_{\alpha\beta}^n e^{2\pi i n\cdot y} ]=
\Lambda_m(n)^{-1}\sum_{|\beta|=m}
(-g_{\beta}^nn^\alpha n^\beta+g_{\alpha}^nn^\beta n^\beta)
=}\atop
\ds{-
\Lambda_m(n)^{-1}n^\alpha\sum_{|\beta|=m}g_{\beta}^n n^\beta
+
\Lambda_m(n)^{-1}g_{\alpha}^n\sum_{|\beta|=m}n^\beta  n^\beta
\stackrel{(\ref{lp5})}=g_{\alpha}^n \Lambda_m(n)^{-1}\sum_{|\beta|=m}n^\beta n^\beta=g_{\alpha}^n,
}
\]
whence the property (\ref{lp4}) follows. The lemma is proved.

\bigskip
We give  here 
another representation lemma which is rather common in homogenization.
\begin{lemma}\label{lem4.2} Let
$g$ be a 1-periodic scalar function from $L^2(Y)$ 
such that $\langle g\rangle=0$.
Then there is a 1-periodic vector $G$ from $H^m(Y)^d$
such that $g=\div \,G$ and
\[
\|G\|_{H^1(Y)}\le c \|g\|_{L^2(Y)}.
\]
\end{lemma}
\doc The required representation is obtained if we set $G=\nabla U$ where $U$ is a solution of the following periodic
problem with laplacian $\Delta=\div\nabla$
\[
U\in H^2(Y),\quad \Delta U=g.
\]
The solvability of this problem is readily shown by using Fourier series. Further details are omitted.

\section{  Discrepancy of the first approximation}
\setcounter{equation}{0}
  \setcounter{pr}{0}
 
 We are aimed to prove the estimate (\ref{13}). It means   
 that the solution $u$ of the homogenized equation 
 approximates  the solution $\ue$ of the initial equation in $L^2$-norm. Therefore,  case 
 the function $u$ is called the zero approximation
 to keep  distinct from the first approximation which approximates the solution $\ue$  in  respect to 
 the Sobolev norm natural to the equation. In our case, this is $H^m$-norm.
 It is quite in common for homogenization theory to use approximations in Sobolev norms to obtain $L^2$-estimate for the difference $\ue{-}u$  as a corollary (see \cite{BLP}, \cite{BP}, \cite{JKO}).
 We shall do the same and try for the first approximation the function
\beq\label{d1}
v^\e(x)=u(x)+\e^m
 \sum\limits_{|\gamma|\le m}N_\gamma^\e(x)D^\gamma u(x),\quad N_\gamma^\e(x)=N_\gamma(y)|_{y=x/\e},
\eeq 
which is a sum of the zero approximation and a corrector term. Here, $N_\gamma(y)$ is a solution to the cell problem  (\ref{c1}). 

To facilitate further actions assume at the first step that the right-hand side function in (\ref{11}) and (\ref{12}) is smooth and with compact support,
that is, $f{\in}\C0(\rd)$. In this case,
 the solution $u$ of the elliptic equation with constant coefficients is smooth and decays exponentially at infinity.
 As a result, $v^\e$ belongs to the space $H^m(\rd)$.
 
 Our goal is to evaluate a discrepancy of $v^\e$ to the equation (\ref{11}). To this end we, first, compare the generalized gradients
 \beq\label{d2}
\Gamma_\a(v^\e,A^\e)=
 \sum\limits_{|\beta|\le  m}a^\e_{\alpha\beta}(y)D^{\beta}v^\e,\quad 
 \Gamma_\a(u,\hat{A})=
 \sum\limits_{|\beta|\le m}\hat{a}_{\alpha\beta} D^{\beta}u,\quad |\a|\le m,
\eeq 
with each other. Easy calculations give
 \beq\label{d3}\ds{
\Gamma_\a(v^\e,A^\e)\stackrel{(\ref{d2}),(\ref{d1})}=
 \sum_{|\beta|\le  m}a^\e_{\alpha\beta}(y)D^{\beta}
( u+
\e^m \sum_{|\gamma|\le m}N_\gamma^\e D^\gamma u)=}\atop\ds{
\sum_{|\beta|\le  m}\left( a^\e_{\alpha\beta}D^{\beta}u+
 \sum_{|\gamma|=m}a^\e_{\alpha\gamma}(\pa^{\gamma}N_\beta)^\e D^{\beta}u\right)+w^\e_\a,}
\eeq 
where $\pa^{\gamma}N_\beta{=}D_y^{\gamma}N_\beta(y)$. The term $w^\e_\a$ collects all the summands, in which there occur derivatives $\pa^{\gamma}N_\beta$ with $|\gamma|{<}m$ and, thereby, there stand 
multipliers $\e^k$, $k{\ge}1$.

Using the notation from (\ref{c5}),(\ref{c6}) we rewrite
\[
\ds{
\Gamma_\a(v^\e,A^\e)=
 \sum_{|\beta|\le  m}\sum_{|\gamma|=m}a^\e_{\alpha\gamma}(e_{\beta\gamma}+(\pa^{\gamma}N_\beta)^\e) D^{\beta}u+w^\e_\a
=
\sum_{|\beta|\le  m}
\tilde{a}^\e_{\alpha\beta}D^{\beta}u+w^\e_\a=}\atop\ds{
\sum_{|\beta|\le  m}
(\tilde{a}^\e_{\alpha\beta}-\hat{a}_{\alpha\beta})D^{\beta}u+\sum_{|\beta|\le  m}
\hat{a}_{\alpha\beta}D^{\beta}u+
w^\e_\a\stackrel{(\ref{d2})}=
 }
 \]
\[
\sum_{|\beta|\le  m}
(\tilde{a}^\e_{\alpha\beta}-\hat{a}_{\alpha\beta})D^{\beta}u+\Gamma_\a(u,\hat{A})+
w^\e_\a.
\]
Finally, 
 \beq\label{d4}\ds{
\Gamma_\a(v^\e,A^\e)=\Gamma_\a(u,\hat{A})+
 \sum_{|\beta|\le  m}g^\e_{\alpha\beta}D^{\beta}
 u+w^\e_\a,
}\atop\ds{
g_{\alpha\beta}(y)=\tilde{a}_{\alpha\beta}(y)-\hat{a}_{\alpha\beta}.
}
\eeq 
Now we transform $g_{\alpha\beta}(y)$
in an appropriate way differently for $\a$ with $|\a|{=}m$ and $|\a|{<}m$.

For fixed $\beta$ the vector $\{g_{\alpha\beta}\}_{|\a|=m}$ satisfies the conditions of Lemma \ref{lem4.1}. Therefore,
there exists the matrix  $\{G_{\alpha\gamma\beta}\}_{|\a|=|\gamma|=m}$ such that
\[g_{\alpha\beta}(y)=\sum_{|\gamma|=m}\pa^\gamma G_{\alpha\gamma\beta}(y),\quad
G_{\alpha\gamma\beta}(y)=-G_{\gamma\alpha\beta}(y), 
\]
and the $H^m$-estimate of the form (\ref{lp3}) is valid for $ G_{\alpha\gamma\beta}(y)$. Hence,
 \beq\label{d5}
 \ds{
g^\e_{\alpha\beta}D^{\beta}u =\sum_{|\gamma|=m} D^\gamma (\e^m G^\e_{\alpha\gamma\beta})D^{\beta}u=}
\atop\ds{
\sum_{|\gamma|=m} D^\gamma (\e^m G^\e_{\alpha\gamma\beta} D^{\beta}u)+w^\e_{\alpha\beta},\quad |\a|{=}m,}
\eeq
where 
$w^\e_{\a\b}$ collects  terms in which there occur derivatives $D^\gamma (\e^m G^\e_{\alpha\gamma\beta})$, $|\gamma|{<}m$, and which, thereby,  contain multipliers $\e^k$, $k{\ge}1$.

As for the coefficients $g^\e_{\alpha\beta}$, $|\a|{<}m$, $|\b|{\le}m$ , from (\ref{d4}), we apply Lemma \ref{lem4.2} to them. So, there exist 1-periodic vectors $\tilde{g}_{\alpha\beta}{\in}H^1_\per(Y)^d$ such that
\[
g_{\alpha\beta}(y)=\div_y \tilde{g}_{\alpha\beta}(y),\quad 
g^\e_{\alpha\beta}(x)=\div_x (\e \tilde{g}^\e_{\alpha\beta}(x)).
\]
Thus,
 \beq\label{d6}
 \ds{
g^\e_{\alpha\beta}D^{\beta}u =\div_x (\e \tilde{g}^\e_{\alpha\beta}(x)) D^{\beta}u=
\div_x (\e \tilde{g}^\e_{\alpha\beta}(x) D^{\beta}u)+\tilde{w}^\e_{\a\b},}
\atop\ds{\tilde{w}^\e_{\a\b}=-\e \tilde{g}^\e_{\alpha\beta}\cdot\nab  D^{\beta}u.
}
\eeq

Since $f=(\hat {A}+\lambda)u$,
we deduce
 \beq\label{d7}
( A^\e+\lambda) v^\e-f=( A^\e+\lambda) v^\e-(\hat {A}+\lambda)u=(A^\e v^\e-\hat {A}u)+\lambda(v^\e-u)
\stackrel{(\ref{d2})}=
 \eeq
 \[
\sum_{|\a|\le m}(-1)^{|\a|} D^\a(\Gamma_\a(v^\e,A^\e)-\Gamma_\a(u,\hat{A}))+\lambda(v^\e-u)
\stackrel{(\ref{d4})} =
 \]
\[
\sum_{|\a|\le m}(-1)^{|\a|} D^\a[
 \sum_{|\beta|\le  m}g^\e_{\alpha\beta}D^{\beta}
 u+w^\e_\a]+\lambda(v^\e-u)=
\]
\[
\sum_{|\a|= m}(-1)^{|\a|} D^\a[\ldots]+\sum_{|\a|< m}(-1)^{|\a|} D^\a[\ldots]+\lambda(v^\e-u)
\stackrel{(\ref{d5}),(\ref{d6})}=
\]
 \beq\label{d8}
{=}(-1)^m \sum_{|\beta|\le  m} \sum_{|\a|=|\gamma|= m}D^\a
D^\gamma (\e^m G^\e_{\alpha\gamma\beta}D^{\beta}u){+}
\sum_{|\beta|\le  m,|\a|< m} 
(-1)^{|\a|} D^\a\div_x (\e \tilde{g}^\e_{\alpha\beta}(x) D^{\beta}u){+}\ldots
 \eeq
 Here, the last dots stand for the terms containing explicit multipliers $\e^k$, $k{\ge}1$, they appear from expressions $v^\e-u$, $w^\e_\a$, $w^\e_{\alpha\beta}$, $\tilde{w}^\e_{\a\b}$ (see (\ref{d1}), (\ref{d4})--(\ref{d6})).
 
 Being outwardly the most complicated, the first sum in (\ref{d8}) presents, actually, the zero functional.  In fact, for any fixed $\b$  we have 
  \beq\label{d9}
F^\e_\b=\sum_{|\a|=|\gamma|= m}D^\a
D^\gamma (\e^m G^\e_{\alpha\gamma\beta}D^{\beta}u)\in H^{-m}(\rd)
 \eeq
 due to the sufficient regularity of the functions $G^\e_{\alpha\gamma\beta}$ and $D^{\beta}u$. Moreover,
  \beq\label{d10}
\langle F^\e_\b,\vp\rangle=\sum_{|\a|=|\gamma|= m}(\e^m G^\e_{\alpha\gamma\beta}D^{\beta}u, D^\a
D^\gamma \vp)=0\quad \forall\, \vp\in \C0(\rd),
 \eeq 
 thanks to the symmetry properties of matrices $\{ G^\e_{\alpha\gamma\beta}\}_{\alpha\gamma}$ and 
 $\{ D^\alpha D^\gamma\vp\}_{\alpha\gamma}$.
 
 Since 
  $f=( A^\e+\lambda) u^\e$  and, thus,
\[
( A^\e+\lambda) v^\e-f=( A^\e+\lambda) v^\e-( A^\e+\lambda) u^\e= A^\e (v^\e-u^\e)+\lambda(v^\e-u^\e),
 \]
 we derive from (\ref{d7})-(\ref{d10}) the equation
  \beq\label{d11}
 \ds{
  A^\e z^\e+\lambda z^\e=\sum_{|\a|\le m}
(-1)^{|\a|} D^\a f^\e_\a
,
}
\atop\ds{z^\e=v^\e-u^\e.
}
\eeq
Here functions $f^\e_\a$
have the 
structure of products
 \beq\label{d12}
\e^kU(x)b(x/\e),\quad k\ge 1, \quad U(x)=D^\gamma u(x),\,|\gamma|\le 2m,\quad 
\eeq
and for 1-periodic functions $b(y)$
there stand  the following expressions
 \beq\label{d13}
\pa^\gamma N_\b(y),\quad \pa^\delta G_{\a\b\gamma}(y),\quad \tilde{g}_{\alpha\beta}(y),\quad |\a|\le m, |\beta|\le m,
 |\gamma|\le m, |\delta|\le m,
\eeq
from the former transformations. By construction, all functions $b(y)$ belong to $L^2(Y)$.
According to the estimate of the form (\ref{6}), the solution of (\ref{d11}) satisfies the inequality
 \beq\label{d14}
\|z^\e\|_{H^m(\rd)}\le C \sum_{|\a|\le m}
\|f^\e_\a\|_{L^2(\rd)}.
\eeq
Here, the majorant is obviously of order $\e$ (see (\ref{d12})), but it
 cannot  be replaced  
 with the expression $\e c_0\|f\|_{L^2(\rd)}$, $c_0{=}const(\la_0, \la_1)$, desired in (\ref{13}).
It would be possible if $b{\in}L^\infty(Y)$, that is not true under our assumptions in general.
In the sequel, we 
show how to overcome this difficulty by introducing an additional parameter of integration.

\section{  Estimate averaged over shifting and its corollaries}
\setcounter{equation}{0}
  \setcounter{theorem}{0}
  \setcounter{lemma}{0} 

\n \textbf{6.1.}
Consider a family of perturbated problems 
  \begin{equation}\label{2.27}
  \begin{array}{cc}
  u_\omega^\e\in  H^m(\rd),  \quad     A_\omega^\e u_\omega^\e+u_\omega^\e=
     f(x),  \quad f\in L^2(\rd),& \\
A_\omega^\e {=}
\sum\limits_{|\alpha|,|\beta|\le m}(-1)^{|\a|}D^{\alpha}(a_{\alpha\beta}(x/\e+\omega)D^{\beta}),  
& \\
  \end{array}
\end{equation}
with shifting parameter $\omega{\in}Y$
in coefficients of the operator. Clearly, (\ref{12})   is 
the corresponding homogenized problem for each $\omega{\in}Y$, the cell problems (see (\ref{c1}))
contain shifting parameter $\omega$ in coefficients, 
and thereby the first approximation for $u_\omega^\e$ is of the form (see (\ref{d1}))
 \begin{equation}\label{2.28}
 v_\omega^\e(x)=u(x)+\e^m
 \sum\limits_{|\gamma|\le m}N_\gamma(y+\omega)D^\gamma u(x),\quad y=x/\e,
\end{equation}
There holds the  corresponding estimate of the form (\ref{d14}) with right-hand side functions defined in  (\ref{d12}), (\ref{d13}). Namely, 
 \[
\|v_\omega^\e-u_\omega^\e\|^2_{H^m(\rd)}\le c\e^2\sum_\a\int_{\rd}|b_\a(\f{x}{\e}+\omega)|^2|U_\a(x)|^2\dx.  
\]
Integrating  in $\omega{\in}Y$ leads to 
\[\ds{
\ilb\|v_\omega^\e-u_\omega^\e\|^2_{H^m(\rd)}d\omega\le c\e^2\sum_\a \ilb\int_{\rd}|b_\a(\f{x}{\e}+\omega)|^2|U_\a(x)|^2\dx d\omega\le}
\atop\ds{c\e^2\sum_\a \|b_\a\|^2_{L^2(Y)}\|U_\a\|^2_{L^2(\rd)}\le c_0\e^2
\|f\|^2_{L^2(\rd)}.
}
\]
Here, at the first step the order of integration is changed, after which we can extract from the integral over $\rd$
 $L^2$-norm of oscillating functions $b_\a$, more exactly, the expression 
 $$\ilb| b_\a(\f{x}{\e}+\omega)|^2d\omega= \|b_\a\|^2_{L^2(Y)}.
  $$ 
  Then we apply the estimate (\ref{ele}) and the estimate
\[
\|b\|_{L^2(Y)}\le c,\, c=const(\lambda_0,\la_1).
\]
The latter is enabled by properties of functions (\ref{d13}). 

Thus, the following lemma is proved.
\begin{lemma}\label{lem6.1} 
Let
$u_\omega^\e$ be a  solution of (\ref{2.27}) and let $v_\omega^\e$ be a function from  (\ref{2.28}). Then 
there holds an integrated (in $\omega{\in}Y$)  
estimate
\begin{equation} \label{2.29}
 \ilb\int_{\rd}(\sum\limits_{1\le|\alpha|\le m}|D^\a(u_\omega^\e-v_\omega^\e)|^2+|u_\omega^\e-v_\omega^\e|^2)\dx d\omega\le c_0\e^2 \|f\|^2_{L^2(\rd)},\quad c_0=const(\lambda_0,\la_1).
\end{equation}
\end{lemma}

It is useful to have 
 another version of integrated (in $\omega{\in}Y$)   
estimate. Take the solution $u^\e(x)$ of the problem (\ref{11}) 
and consider a family of shifted functions
$\tilde{u}_\omega^\e(x){=}u^\e(x{+}\e\omega)$.
They satisfy the equation (\ref{2.27}) with the shifted right-hand side function $f(x{+}\e\omega)$, 
\[
 A_\omega^\e\tilde{u}_\omega^\e+\tilde{u}_\omega^\e=f(x{+}\e\omega).
\]
Then
\begin{equation} \label{2.30}
 A_\omega^\e(u_\omega^\e-\tilde{u}_\omega^\e)+u_\omega^\e-\tilde{u}_\omega^\e=f(x)-f(x{+}\e\omega).
\end{equation}
By properties of shifting, 
$$
(\ild(f(x)-f(x{+}\e\omega))\varphi(x)\dx
)^2\le \|f\|^2_{L^2(\rd)}\ild|\varphi(x)-\varphi(x{-}\e\omega)|^2\dx\le
$$
$$
c\e^2\|f\|^2_{L^2(\rd)}\|\nab \varphi\|^2_{L^2(\rd)},\quad c=const(d),
$$
and therefore,  from (\ref{2.30}) it readily follows that  
\[
\|u_\omega^\e-\tilde{u}_\omega^\e\|_{H^m(\rd)}\le \e C
\|f\|_{L^2(\rd)}.
\]
Thus,   $u_\omega^\e$ may be replaced
with $\tilde{u}_\omega^\e$ 
in (\ref{2.29}) without detriment to the right-hand side of (\ref{2.29}).
Namely,
\begin{equation} \label{2.31}
 J{:=}\ilb\int_{\rd}
 \sum\limits_{1\le|\alpha|\le m}|D^\a(\tilde{u}_\omega^\e-v_\omega^\e)|^2\dx d\omega{+}\ilb\int_{\rd}|\tilde{u}_\omega^\e-v_\omega^\e|^2\dx d\omega{\le} c\e^2 \|f\|^2_{L^2(\rd)}, 
\end{equation} 
where $c=const(\lambda_0,\la_1).$

\n \textbf{6.2.}  Now we are going to derive some corollaries from (\ref{2.31}). 

1$^\circ$ 
Discarding the first integral in
 (\ref{2.31}) and changing the order of integration in the remaining one, we deduce, by convexity,
 \begin{equation} \label{2.32}
 \int_{\rd}|
 \langle\tilde{u}_\omega^\e-v_\omega^\e\rangle_\omega|^2\dx \le c\e^2 \|f\|^2_{L^2(\rd)},
\end{equation}
where $\langle\cdot\rangle_\omega{=}\ilb\cdot\,d\omega$.
Clearly,
\[
\langle v_\omega^\e\rangle_\omega\stackrel{(\ref{2.28})}=u(x),\quad
\langle \tilde{u}_\omega^\e\rangle_\omega=\ilb u^\e(x+\e\omega)\,d\omega=(S^\e u^\e)(x)
\]
is Steklov average of the function $u^\e(x)$. We recall the following property
\begin{equation} \label{2.33}
\|S^\e\varphi-\varphi\|_{L^2(\rd)} \le c\e \|\nab\varphi\|_{L^2(\rd)},\quad c=const(d),
\end{equation}
for the Steklov average of the function $\varphi$ defined as
\[
(S^\e\varphi)(x)=\ilb \varphi(x+\e\omega)\,d\omega.
\]
Therefore, 
(\ref{2.32}) means
\begin{equation} \label{2.340}
\|S^\e u^\e-u\|_{L^2(\rd)} \le c\e \|f\|_{L^2(\rd)},
\end{equation}
and, by triangle inequality,
\[
\|u^\e-u\|_{L^2(\rd)}\le \|u^\e-S^\e u^\e\|_{L^2(\rd)}+
\|S^\e u^\e-u\|_{L^2(\rd)} \le c_0\e \|f\|_{L^2(\rd)}.
\]
Here the property (\ref{2.33}) of Steklov averaging is applied to
  $u^\e$ and, finally, the evident inequality  $\|\nab u^\e\|_{L^2(\rd)}\le C\|f\|_{L^2(\rd)}$, arising from the energy estimate,
is used. 
As a result, the estimate
(\ref{9}) is proved.

From (\ref{2.29}), 
we derive  similarly 
the inequality 
\[
\sum\limits_{1\le|\alpha|\le m}\|D^\a S^\e u^\e-D^\a u\|_{L^2(\rd)} \le c\e \|f\|_{L^2(\rd)},
\]
which together with (\ref{2.340}) leads to 
\begin{equation} \label{2.341}
\|S^\e u^\e-u\|_{H^m(\rd)} \le c_0\e \|f\|_{L^2(\rd)},\quad c_0 =const(\la_0, \la_1).
\end{equation}
This $H^2$-estimate deserves attention 
because it does not contain  any corrector. The 
proximity between  $u^\e$ and $u$ in $H^2$-norm is achieved via Steklov smoothing only.

2$^\circ$ Now transform the expression $J$ from (\ref{2.31}) in another fashion. First, change the variable
$x\to x^\prime{=}x{+}\e\omega$,
which leads to 
$$\tilde{u}_\omega^\e(x){=}u^\e(x^\prime),$$
$$\tilde{u}_\omega^\e(x)-v_\omega^\e(x)=
u^\e(x^\prime)-
u(x^\prime-\e\omega)-\e^m
 \sum\limits_{|\gamma|\le m}N_\gamma(x^\prime/\e)D^\gamma u(x^\prime-\e\omega)
, $$
$$
J\stackrel{(\ref{2.31})}=
\ilb\int_{\rd}[|{u}^\e(x){-}u(x{-}\e\omega){-}
\e^m
 \sum\limits_{|\gamma|\le m}N_\gamma(x/\e)D^\gamma u(x-\e\omega)
|^2{+}|\sum\limits_{1\le|\alpha|\le m}|D^\a
(\ldots)|^2]\dx d\omega,
$$
where the dots  keep back the difference of functions from the previous summand.

Hence, after changing the order of integration we deduce, by convexity,
that
\begin{equation} \label{2.34}
J\ge \int_{\rd}[|z^\e|^2{+}\sum\limits_{1\le|\alpha|\le m}|D^\a z^\e|^2]\dx,
\end{equation}
where 
\begin{equation} \label{2.35}
\ds{z^\e(x)=
u^\e(x)-\langle u(x{-}\e\omega)\rangle_\omega-\e^m
 \sum\limits_{|\gamma|\le m}N_\gamma(x/\e)\langle D^\gamma u(x-\e\omega)\rangle_\omega=
}\atop\ds{
u^\e(x)-(S^\e u)(x)-\e^m
 \sum\limits_{|\gamma|\le m}N_\gamma(x/\e)S^\e( D^\gamma u)(x).
}
\end{equation}
So, (\ref{2.31}) and (\ref{2.34}) imply the estimate
\[
\int_{\rd}[|z^\e|^2{+}\sum\limits_{1\le|\alpha|\le m}|D^\a z^\e|^2]\dx\le c\e^2\|f\|_{L^2(\rd)}^2.
\]
Here, without detriment to the right-hand side, one can replace 
in $z^\e$
(see (\ref{2.35})) the  Steklov average $S^\e u$
with the function $u$ itself,
having in mind the property (\ref{2.33}) and the elliptic estimate for $u$.
In such a way there appear the first approximation with smoothed corrector (smoothig in Steklov sense)
 \begin{equation} \label{2.36}
\hat{v}^\e(x)= u(x)+\e^m
 \sum\limits_{|\gamma|\le m}N_\gamma(x/\e)S^\e( D^\gamma u)(x),
\end{equation}
and the estimate
\[
\int_{\rd}[|u^\e-\hat{v}^\e|^2{+}
\sum\limits_{1\le|\alpha|\le m}|D^\a
(u^\e-\hat{v}^\e))|^2]
\dx\le c\e^2\|f\|_{L^2(\rd)}^2.
\]

As a result, the following theorem about approximation in $H^m$-norm is proved. 
\begin{theorem}\label{th6.2}
For the difference of the solution $u^\e$ to the problem (\ref{11})  and the function $\hat{v}^\e$, defined
in (\ref{2.36}),
there holds the estimate
\begin{equation} \label{2.70}
 \|u^\e-\hat{v}^\e\|_{H^m(\rd)}\le c_0\e \|f\|_{L^2(\rd)},\quad c_0 =const(\la_0, \la_1).
\end{equation}
\end{theorem}

Accordingly,  (\ref{2.36}) and (\ref{2.70}) imply the operator-type estimate (\ref{14}) with the correcting operator
$$\mathcal{K}^\e=\e^m
 \sum\limits_{|\gamma|\le m}N_\gamma(x/\e)S^\e D^\gamma(\hat{A}+\la)^{-1}.
$$

\section{  Some remarks}
\setcounter{equation}{0}
  \setcounter{theorem}{0}
  \setcounter{lemma}{0} 

\textbf{1.}
In our method, properties of shifting and Steklov averaging (or smoothing) are essential. We omit here their proof, all the necessary proofs are given 
  in
\cite{Z05},\cite{ZP}. In particular, by properties of Steklov average and due to the special structure, the function $\hat{v}^\e$, defined
in (\ref{2.36}), belongs to the Sobolev space $H^m(\rd)$ which is, of course, necessary for the estimate (\ref{2.70}).
Note, that we have gained this property of $\hat{v}^\e$ automatically as a byproduct while deriving the estimate (\ref{2.70}).
 For all that, we have also used some  properties of Steklov  smoothing.
We can prove directly that   $\hat{v}^\e{\in}H^m(\rd)$. 
For, this function itself and all its derivatives up to the order $m$ are elements of $H^m(\rd)$ due to following assertion proved in \cite{ZP}:
\textit{suppose that  $b{\in}L^2(Y)$, 
 $b$ is 1-periodic, $b_\e(x){=}b(x/\varepsilon)$ and $\varphi{\in}L^2(\rd)$, 
 then \\
 $\quad{}\quad$ $\quad{}\quad$  $\quad{}\quad$  $\quad{}\quad$ $\quad{}\quad$ 
 $\|b_\e S^\varepsilon \varphi\|_ {L^2(\rd)}\le \|b\|_{L^2(Y)}\|\varphi\|_ {L^2(\rd)}$.
 } 
 
 In general, under our minimal assumptions, when coefficients  of the operator 
 are only measurable bounded functions and the right-hand side function $f$ belongs to $L^2(\rd)$,  the definition (\ref{2.36}) without Steklov smoothing in it (that is exactly  (\ref{d1}))
does not enable $H^m$-regularity  of the approximation, 
at least, from the first sight. There are some particular cases when Steklov smoothing can be omitted in (\ref{2.36}) and, thus, the estimate  (\ref{2.70}) is also true with the approximation ${v}^\e$ from (\ref{d1}) instead of $\hat{v}^\e$. The examples are given below without going into details, for, the full justification may be not obvious and even
cumbersome. Instead of detailing, we make  reference to our papers, if possible.

\bigskip
\n
{E x a m p l e 1} (general problem in dimension $d{=}1$). In one-dimensional case, the cell problems are solved explicitly, the cell functions, with all their derivatives up to order $m$, are bounded; thereby  the  justification is easy.

\n
E x a m p l e 2 (general problem for the order $2m{=}2$). Second-order operators of the form (\ref{8}) can be treated by arguments  considered  for more particular case in \cite{ZP}.

\n
E x a m p l e 3 (operator with bilaplacian). Fourth-order operators of the form $A^\e=\Delta a(x/\e)\Delta$, where $\Delta$  is $d$-dimensional Laplacian and $a(y)$ is a 
positive  function from $L^\infty_\per(Y)$, produce a very peculiar cell problem 
which leads to rather simple  in form the first approximation. This case is  considered in \cite{P15}.

\textbf{2.} Now we give examples of operators satisfying the condition (\ref{3}).
 Consider an operator
\beq\label{3.600}
A=\sum\limits_{i,j,s,h}\f{\partial^2}{\partial x_i\partial x_j}(a_{ijsh}(x)
\f{\partial^2}{\partial x_s\partial x_h})
\eeq
with a fourth order tensor  $a(x){=}\{a_{ijsh}(x)\}$ acting in a space of $(d{\times}d)$-matrices. Assume that
\beq\label{3.601}
a_{ijsh}{=}a_{shij},\, a_{ijsh}{=}a_{jish}{=} a_{ijhs}.
\eeq
So, the tensor $a$
defines a symmetric operator in the space of  symmetric  matrices.
\\
\,\,C\,a\,s\,e 1. Assume also that 
\beq\label{3.60}
\lambda_0\xi\cdot\xi\le a\xi\cdot\xi\le\lambda_0^{-1}\xi\cdot\xi 
\eeq
 with a constant $\lambda_0{>}0$ for any symmetric  matrix $\xi{=}\{\xi_{ij}\}$, where $\xi\cdot\eta{=}\xi_{ij}\eta_{ij}$
(here and hereafter, over repeated indices summation is assumed from 1 to d). 
The left inequality (\ref{3.60}) enables (\ref{3}), since for 
$ u\in \C0(\rd)$  the point-wise inequality
\[
a_{ijsh}\f{\partial^2 u}{\partial x_s\partial x_h}\f{\partial^2 u}{\partial x_i\partial x_j}\ge
\lambda_0 \f{\partial^2 u}{\partial x_i\partial x_j} \f{\partial^2 u}{\partial x_i\partial x_j}
\]
is valid which after integration  leads to
 (\ref{3}). Operators of the type (\ref{3.600}) satisfying the conditions (\ref{3.601}) and  
 (\ref{3.60}) appear in the elasticity theory  applied to thin plates.
 \\
\,\,C\,a\,s\,e 2.
Consider a tensor $a$ acting on a matrix $\xi$
 as $a\xi{=}\alpha(Tr\,\xi)E$, where $Tr\,\xi=\xi_{ii}$, $E$ is an identity matrix, $\alpha{\in}L^\infty(\rd)$, $\alpha{\ge}\alpha_0{>}0$. 
Obviously, $a\xi\cdot\xi{=}\alpha(Tr\,\xi)^2$. Thus,
for the matrix $\xi{=}D u{=}\{\f{\partial^2 u}{\partial x_i\partial x_j}\}$, we have  
$Tr\,D u{=}\Delta u$, $aD u\cdot D u{=}\alpha \Delta u \Delta u.$
The corresponding operator is $A{=} \Delta  (\alpha(x) \Delta )$. Note that the left inequality in (\ref{3.60}) does not hold for the tensor $a$, though the estimate (\ref{3}) is true:
\[
\ild \alpha(x) \Delta u \Delta u\,dx\ge \alpha_0\ild |\Delta u|^2\,dx\ge \lambda_0\ild \sum_{i,j}
|\f{\partial^2 u}{\partial x_i\partial x_j}|^2\dx.
\]
The last inequality with the appropriate  constant $ \lambda_0{>}0$  is verified via Fourier transformation. 

\textbf{3.}   
After the paper with these results  was submitted to  "Applicable Analysis" in September 2015 ${}$there appeared the publication  on the close topic 
 \cite{KS} related to matrix strongly elliptic self-adjoint operators.

\textbf{4.}
The method used here to prove the operator-type estimates in periodic homogenization
allows to extend our results to elliptic operators  with locally periodic or multiscale coefficients
 (see the appropriate technique in \cite{MIAN}, \cite{FA}).
  As for the spectral approach used in \cite{V} and \cite{KS}, this is impossible.
We can treat also elliptic operators with complex-valued coefficients. Of course, 
the condition  
of ellipticity  (\ref{3}) need to 
be slightly modified in this case. Here, the  case of real-valued  coefficients is chosen only for simplicity.

\bigskip 
\noindent\textbf{ Acknowledgements}.
The author was supported by
 Russian Science Foundation (grant no. 14-11-00398).

\end{document}